%% file: tra4.tex
\documentclass[11pt]{article}
\usepackage{psfig}
\usepackage{latexsym}
\usepackage{amsmath}
\usepackage[usenames]{color}
\usepackage[dvips]{epsfig}
\usepackage{graphicx}

\newcommand{\Cor}{\mbox{Cor}}
\newcommand{\COR}{\mbox{COR}}

\newcommand{\mychoose}[2]{ \left(\!\!\begin{array}{c}
   #1 \\ #2 \end{array}\!\!\right)}

\newcommand{\itb}{\begin{itemize}}
\newcommand{\itn}{\end{itemize}}
\newcommand{\nmn}{\end{enumerate}}
\newcommand{\nmb}{\begin{enumerate}}
\newcommand{\ctb}{\begin{center}}
\newcommand{\ctn}{\end{center}}
\newcommand{\tb}{\begin{center}\Large \bf}
\newcommand{\tn}{\end{center}}
\newcommand{\arb}{\begin{array}}
\newcommand{\qrn}{\end{array}}
\newcommand{\tblb}{\begin{tabular}}
\newcommand{\tbln}{\end{tabular}}
\newcommand{\fb}{\begin{figure}}
\newcommand{\fn}{\end{figure}}

\def\bern{\mbox{Bernoulli}\,}
\def\betad{\mbox{Beta}\,}

\def\ybar{\overline{y}}

\def\beq{\begin{eqnarray}}
\def\eeq{\end{eqnarray}}
\def\mb#1{\mbox{\boldmath $#1$}}

\def\x{\mb x}

\def\thetab{\mb \theta}
\def\phib{\mb \phi}

\def\IID{
 \,\begin{array}{cc} \\[-20pt] \mbox{\tiny IID} \\[-7pt] \sim \end{array}\,}
\def\trn{^{\mbox{\tiny train}}}

\begin{document}

\fontsize{11}{14.5pt}\selectfont

\begin{center}

{\small Technical Report No.\ 0705,
 Department of Statistics, University of Toronto}

\vspace*{0.95in}

{\Large\bf A Method for Avoiding Bias from Feature Selection with\\[6pt]
           Application to Naive Bayes Classification Models}
\\[15pt]

{\large Longhai Li, Jianguo Zhang, and Radford M. Neal}\\[3pt]
 Department of Statistics \\
 University of Toronto, Toronto, Ontario, Canada \\
 \texttt{\{longhai,jianguo,radford\}@stat.utoronto.ca}\\[10pt]

 19 February 2007
\end{center}

\vspace{8pt}

\noindent \textbf{Abstract.}  For many classification and regression
problems, a large number of features are available for possible use
--- this is typical of DNA microarray data on gene expression, for
example.  Often, for computational or other reasons, only a small
subset of these features are selected for use in a model, based on
some simple measure such as correlation with the response variable.
This procedure may introduce an optimistic bias, however, in which the
response variable appears to be more predictable than it actually is,
because the high correlation of the selected features with the
response may be partly or wholely due to chance.  We show how this
bias can be avoided when using a Bayesian model for the joint
distribution of features and response.  The crucial insight is that
even if we forget the exact values of the unselected features, we
should retain, and condition on, the knowledge that their correlation
with the response was too small for them to be selected. In this paper
we describe how this idea can be implemented for ``naive Bayes''
models of binary data.  Experiments with simulated data confirm that
this method avoids bias due to feature selection.  We also apply the
naive Bayes model to subsets of data relating gene expression to colon
cancer, and find that correcting for bias from feature selection does
improve predictive performance.

\section{\hspace*{-7pt}Introduction}\label{sec-intro}\vspace*{-8pt}

Regression and classification problems that have a large number of
available ``features'' (also known as ``inputs'', ``covariates'', or
``predictor variables'') are becoming increasingly common.  Such
problems arise in many application areas.  Data on the expression
levels of tens of thousands of genes can now be obtained using DNA
microarrays, and used for tasks such as classifying tumors.  Document
analysis may be based on counts of how often each word in a large
dictionary occurs in each document.  Commercial databases may contain
hundreds of features describing each customer.

Using all the features available in such problems is often infeasible.
Using too many features can result in ``overfitting'' when simple
statistical methods such as maximum likelihood are used, with the
consequence that poor predictions are made for the response variable
(e.g., the class) in new items.  More sophisticated Bayesian methods can
avoid such statistical problems, but using a large number of features
may still be undesirable.  We will focus primarily on situations where
the computational cost of looking at all features is too burdensome.
Another issue in some applications is that using a model that looks at
all features will require measuring all these features when making
predictions for future items, which may sometimes be costly.  In some
situations, models using few features may be preferred because they
are easier to interpret.

For such reasons, modelers often use only a subset of features,
chosen by some simple indicator of how useful they might be in
predicting the response variable --- see, for example, the papers in
(Guyon, \textit{et al.} 2006).  For both regression problems with a
real-valued response variable and classification problems with a
binary (0/1) class variable, one suitable measure of how useful a
feature may be is the sample correlation of the feature with the
response.  If the absolute value of this sample correlation is small,
we might decide to omit the feature from our model.  This criterion is
not perfect, of course --- it may result in a relevant feature being
ignored if its relationship with the response is non-linear, and it
may result in many redundant features being retained even when they
all contain essentially the same information.  Sample correlation is
easily computed, however, and hence is an attractive criterion for
screening a large number of features.

Unfortunately, a model that uses only a subset of features, selected
based on their high correlation with the response, will be
optimistically biased --- i.e., predictions made using the model will
(on average) be more confident than is actually warranted.  For
example, we might find that the model predicts that certain items
belong to class~1 with probability 90\%, when in fact only 70\% of
these items are in class~1.  In a situation where the class is
actually completely unpredictable from the features, a model using a
subset of features that purely by chance had high sample correlation
with the class may produce highly confidence predictions that have
less actual chance of being correct than just guessing the most common class.

This optimistic bias comes from ignoring a basic principle of Bayesian
inference --- that we should base our conclusions on probabilities
that are conditional on \textit{all} the available information.  If we
have an appropriate model, this principle would lead us to use all the
features.  This would produce the best possible predictive
performance.  However, we assume here that computational or other
pragmatic issues make using all features unattractive.  When we
therefore choose to ``forget'' some features, we can nevertheless
still retain the information about how we selected the subset of
features that we use in the model.  Properly conditioning on this
information when forming the posterior distribution eliminates the
bias from feature selection, producing predictions that are as good as
possible given the information in the selected features, without the
overconfidence that comes from ignoring the feature selection process.

In the next section, we describe this idea in more detail, and discuss
the difficulties of implementing it.  We then show how the idea can be
applied to a simple ``naive Bayes'' classification model with binary
features that are assumed to be independent given the value of the
response variable.  We apply this naive Bayes model to simulated data
and to data regarding gene expression in colon cancer, showing that
bias correction does indeed improve predictions.  Our method is more
generally applicable, however.  In the final section, we briefly
discuss our work on mixture models for binary data and on factor
analysis models for real-valued data, as well as other possible
applications.

\section{\hspace*{-7pt}Our method for avoiding selection
                       bias}\label{sec-idea}\vspace*{-8pt}

Suppose we wish to predict a response variable, $y$, based on the
information in the numerical features $x_1,\ldots,x_p$, which we
sometimes write as a vector, $\x$.  Our method is applicable both when
$y$ is a binary ($0/1$) class indicator, as is the case for the naive
Bayes models discussed later, and when $y$ is real-valued.  We assume
that we have complete data on $n$ ``training'' cases, for which the
responses are $y^{(1)},\ldots,y^{(n)}$ (collectively written as
$y\trn$) and the feature vectors are $\x^{(1)},\ldots,\x^{(n)}$
(collectively written as $\x\trn$).  (Note that when $y$, $\x$, or
$x_t$ are used without a superscript, they will refer to some
unspecified case.)  We wish to predict the response for one or more
``test'' cases, for which we know only the feature vector.  Our
predictions will take the form of a distribution for $y$, rather than
just a single-valued guess.

We are interested in problems where the number of features, $p$, is
quite big --- perhaps as large as ten or a hundred thousand --- and
accordingly (for pragmatic reasons) we intend to select a subset of
features based on the absolute value of each feature's sample
correlation with the response.  The sample correlation of the
response with feature $t$ is defined as follows (or as zero if
the denominator below is zero):\vspace*{-13pt}
\beq
  \COR (y\trn,\,x\trn_t) & = &
  {\displaystyle\sum_{i=1}^n\, \big(y^{(i)}-\bar y\big)\,
                                 \big(x^{(i)}_t-\bar x_t\big)
  \over
   \sqrt{\sum\limits_{i=1}^n \big(y^{(i)}-\bar y\big)^2}\
   \sqrt{\sum\limits_{i=1}^n \big(x^{(i)}_t-\bar x_t\big)^2} }
\label{eq-cor-def}
\eeq
where $\bar y\, =\, {1 \over n}\sum\limits_{i=1}^n y^{(i)}$ and $\bar x_t\, =\,
{1 \over n}\sum\limits_{i=1}^n x^{(i)}_t$.  The numerator
can be simplified to $\sum\limits_{i=1}^n \big(y^{(i)}-\bar y\big)
x^{(i)}_t$.

Although our interest is only in predicting the response, we assume
that we have a model for the joint distribution of the response
together with all the features.  From such a joint distribution, with
probability or density function $P(y,x_1,\ldots,x_p)$, we can obtain
the conditional distribution for $y$ given any subset of features, for
instance $P(y\,|\,x_1,\ldots,x_k)$, with $k<p$.  This is the
distribution we need in order to make predictions based on this
subset.  Note that selecting a subset of features makes sense only
when the omitted features can be regarded as random, with some
well-defined distribution given the features that are retained, since
such a distribution is essential for these predictions be meaningful.
This can be seen from the following expression:
\beq
  \lefteqn{P(y\,|\,x_1,\ldots,x_k)}\ \ \ \ \ \ \nonumber\\[5pt]
  & = &
  \int \cdots \int P(y\,|\,x_1,\ldots,x_k,x_{k+1},\ldots,x_p)\
                   P(x_{k+1},\ldots,x_p\,|\,x_1,\ldots,x_k)\
                   dx_{k+1} \cdots dx_p\ \ \ \
\eeq
If $P(x_{k+1},\ldots,x_p\,|\,x_1,\ldots,x_k)$ does not exist in any
meaningful sense --- as would be the case, for example, if the data
were collected by an experimenter who just decided arbitrarily what to
set $x_{k+1},\ldots,x_p$ to --- then $P(y\,|\,x_1,\ldots,x_k)$ will
also have no meaning.

Consequently, features that cannot usefully be regarded as random
should always be retained.  Our general method can accommodate such
features, provided we use a model for the joint distribution of the
response together with the random features, conditional on given
values for the non-random features.  However, for simplicity, we will
ignore the possible presence of non-random features in this paper.

We will assume that a subset of features is selected by fixing a
threshold, $\gamma$, for the absolute value of the correlation of a
selected feature with the response.  We then omit feature $t$ from the
feature subset if $|\COR (y\trn,x\trn_t)|\,\le\,\gamma$, retaining
those features with a greater degree of correlation.  Another possible
procedure is to fix the number of features, $k$, that we wish to
retain, and then choose the $k$ features whose correlation with the
response is greatest in absolute value, breaking any tie at random.
If $s$ is the retained feature with the weakest correlation with the
response, we can set $\gamma$ to $|\COR (y\trn,x\trn_s)|$, and we will
again know that if $t$ is any omitted feature, $|\COR
(y\trn,x\trn_t)|\,\le\,\gamma$.  If either the response or the
features have continuous distributions, exact equality of sample
correlations will have probability zero, and consequently this
situation can be treated as equivalent to one in which we fixed
$\gamma$ rather than $k$.  If sample correlations for different
features can be exactly equal, we should theoretically make use of the
information that any possible tie was broken the way that it was, but
ignoring this subtlety is unlikely to have any practical effect, since
ties are still likely to be rare.

Regardless of the exact procedure used to select features, we will
denote the number of features retained by $k$, we will renumber the
features so that the subset of retained features is $x_1,\ldots,x_k$,
and we will assume we know that $|\COR (y\trn,x\trn_t)|\,\le\,\gamma$
for $t=k\!+\!1,\ldots,p$.

We can now state the basic principle behind our bias-avoidance
method:\ \ When forming the posterior distribution for parameters of
the model using a subset of features, we should condition not only on
the values in the training set of the response and of the $k$ features
we retained, but also on the fact that the other $p\!-\!k$ features
have sample correlation with the response that is less than $\gamma$ 
in absolute value.  That
is, the posterior distribution should be conditional on the following
information:\vspace*{-4pt}
\beq
   y\trn,\ \ \x\trn_{1:k},\ \ |\COR (y\trn,x\trn_t)|\,\le\,\gamma\
                              \mbox{for $t = k\!+\!1,\ldots,p$}
\label{eq-cond-info}
\eeq
where $\x\trn_{1:k}\, =\, (x\trn_1,\ldots,x\trn_k)$.

We claim that this procedure of conditioning on the fact that
selection occurred will eliminate the bias from feature selection.
Here, ``bias'' does not refer to estimates for model parameters, but
rather to our estimate of how well we can predict responses in test
cases.  Bias in this respect is also referred to as a lack of
``calibration'' --- that is, the predictive probabilities do not
represent the actual chances of events (Dawid 1982).  If the model
describes the actual data generation mechanism, and the actual values
of the model parameters are indeed randomly chosen according to our
prior, Bayesian inference always produces well-calibrated results, on
average (with respect to the distribution of data and parameter values
chosen from the prior).

In justifying our claim that this procedure avoids selection bias, we
will assume that our model for the joint distribution of the response
and all features, and the prior we chose for it, are appropriate for
the problem, and that we would therefore not see bias if we predicted
the response using all the features.  Now, imagine that rather than
selecting a subset of features ourselves, after seeing all the data,
we instead set up an automatic mechanism to do so, providing it with
the value of $\gamma$ to use as a threshold.  This mechanism, which
has access to all the data, will compute the sample correlations of
all the features with the response, select the subset of features by
comparing these sample correlations with $\gamma$, and then erase the
values of the omitted features, delivering to us only the identities
of the selected features and their values in the training cases.  If
we now condition on all the information that \textit{we} know, but not
on the information that was available to the selection mechanism but
not to us, we will obtain unbiased inferences.  The information we
know is just that of (\ref{eq-cond-info}) above.

Our method requires computation of an adjustment factor, $P({\cal
S}\,|\,\alpha,\,y\trn)$, where $\alpha$ is the set of parameters whose
likelihood needs adjusting, and ${\cal S}$ represents the information
regarding selection, namely that $|\COR (y\trn,x\trn_t)|\,\le\,\gamma$
for $t = k\!+\!1,\ldots,p$.  Computing this factor is much easier if
the $x\trn_t$ are conditionally independent given $\alpha$ and $y\trn$, since
we can then write it as a product of factors pertaining to the
various omitted features.  For the models we consider, these factors
are also \textit{all the same}, since nothing distinguishes one
omitted feature from another.  We can then write
\beq
    P({\cal S}\,|\,\alpha,\,y\trn) & = &
    \!\!\prod_{t=k+1}^p\!\! P\big(|\COR (y\trn,x\trn_t)|\,\le\,\gamma
                              \,|\,\alpha,\,y\trn) \\[5pt]
    & = &
    \Big[ P\big(|\COR (y\trn,x\trn_t)|\,\le\,\gamma
                              \,|\,\alpha,\,y\trn) \Big]^{p-k}
\label{eq-cor-lik}\eeq
where in the second expression, $t$ represents \textit{any} of the
omitted features.  Note that in this expression, $y\trn$ is
conditioned on, and hence considered fixed, whereas $x\trn_t$ is
random.  Since the time needed to compute this adjustment factor does
not depend on the number of omitted features, we may hope to save a
large amount of computation time by omitting many features.

Computing the single factor we do need is not trivial, however, since
it involves integrals over both $x\trn_t$ and any parameters specific
to particular features.  As we will see, however, efficient
computation is possible for the naive Bayes model.

\section{\hspace*{-7pt}Application to naive Bayes
         models with binary features}\vspace*{-8pt}\label{sec-bnaive}

In this section we show how to apply the bias correction method to
Bayesian naive Bayes models in which both the features and the
response are binary.  Binary features are natural for some problems
(e.g., test answers that are either correct or incorrect), or may result
from thresholding real-valued features.  Such thresholding can
sometimes be beneficial --- in a document classification problem, for
example, whether or not a word is used at all may be more relevant to
the class of the document than how many times it is used.  Naive Bayes
models assume that features are independent given the response.  This
assumption is often incorrect, but such simple naive Bayes models have
nevertheless been found to work well for many practical problems. Here
we show how to correct for selection bias in binary naive Bayes
models, whose simplicity allows the required adjustment factor to be
computed very quickly.  Simulations reported in the next section
show that substantial bias can be present with the
uncorrected method, and that it is indeed corrected by conditioning on
the fact that feature selection occurred.  We then apply the method to
real data on gene expression relating to colon cancer, and again find
that our bias correction method improves predictions.

\subsection{\hspace*{-4pt}Definition of the binary naive Bayes
                          model}\vspace*{-4pt}\label{sub-model-b}

Let $\x^{(i)}=(x^{(i)}_1,\cdots,x^{(i)}_p)$ be the vector of $p$
binary features for case $i$, and let $y^{(i)}$ be the binary response
for case $i$, indicating the class.  For example, $y^{(i)}=1$ might
indicates that cancer is present for patient $i$, and $y^{(i)}=0$
indicate that cancer is not present.  Cases are assumed to be
independent given the values of the model parameters (ie, exchangeable
\textit{a priori}).  The probability that $y=1$ in a case is given by the
parameter $\psi$.  Conditional on the class $y$ in some case (and on
the model parameters), the features $x_1,\ldots,x_p$ are assumed to be
independent, and to have Bernoulli distributions with parameters
$\phi_{y,1},\ldots,\phi_{y,p}$, collectively written as $\phib_y$,
with $\phib=(\phib_0,\phib_1)$ representing all such parameters.
In other words, the data is modeled as
\beq
   y^{(i)}\ |\ \psi & \sim & \bern(\psi),\ \ \ \mbox{for $i=1,\ldots,n$}
   \label{sample-y-b}
\\[4pt]
   x^{(i)}_j\ |\ y^{(i)},\, \phib & \sim & \bern(\phi_{y^{(i)},j}),
   \ \ \ \mbox{for $i=1,\ldots,n$ and $j=1,\ldots,p$}
   \label{sample-x-b}
\eeq

We use a hierarchical prior that expresses the
possibility that some features may have almost the same distribution
in the two classes.  In detail, the prior has the following form:
\beq
 \psi & \sim & \betad(f_1,f_0) \label{prior-psi-b} \\[4pt]
 \alpha & \sim & \mbox{Inverse-Gamma}(a,b)\label{prior-alpha-b}
 \\[4pt]
 \theta_1,\ldots,\theta_p &\IID & \mbox{Uniform}(0,1)
   \label{prior-theta-b}
 \\[4pt]
 \phi_{0,j},\,\phi_{1,j}\ |\ \alpha,\,\theta_j & \IID &
   \betad(\alpha\theta_j,\,\alpha(1\!-\!\theta_j)),
   \ \ \ \mbox{for $j=1,\ldots,p$}
 \label{prior-phi-b}
\eeq
The hyperparameters $\thetab\, =\, (\theta_1,\ldots,\theta_p)$
are used to introduce dependence between $\phib_{0,j}$ and $\phib_{1,j}$,
with $\alpha$ controlling the degree of dependence.
Features for which $\phi_{0,j}$ and $\phi_{1,j}$ differ greatly
are more relevant to predicting the response.  When $\alpha$
is small, the variance of the Beta distribution
in~(\ref{prior-phi-b}), which is 
$\theta_j\,(1\!-\!\theta_j)\,/\,(\alpha\!+\!1)$, is large, and
many features are likely to have predictive power, whereas when $\alpha$ is
large, it is likely that most features will be of little use in predicting 
the response, since $\phi_{0,j}$ and $\phi_{1,j}$ are likely to be 
almost equal.
We chose an Inverse-Gamma prior for $\alpha$ (with density function
proportional to $\alpha^{-(1+a)}\exp(-b/\alpha)$) because it has
a heavy upward tail, allowing for the possibility that $\alpha$ is 
large.  Our method of correcting selection bias will have the
effect of modifying the likelihood in a way that favors larger values
for $\alpha$ than would result from ignoring the effect of selection.

\subsection{\hspace*{-4pt}Integrating away $\psi$ and 
            $\phib$}\vspace*{-4pt}\label{sub-int-b}

Although the above model is defined with $\psi$ and $\phib$ parameters for
better conceptual understanding, computations are simplified by integrating
them way analytically.  

Integrating away $\psi$, the joint probability of
$y\trn = (y^{(1)},\ldots,y^{(n)})$ is as follows, 
where $I(\,\cdot\,)$ is the indicator function, equal to 1 if the enclosed 
condition is true and 0 if it is false:
\beq
  P(y\trn)
  & = & \int_0^1 {{\Gamma(f_0+f_1)}\over{\Gamma(f_0)\Gamma(f_1)}}
                \psi^{f_1}(1\,-\,\psi)^{f_0}\
        \psi^{\sum\limits_{i=1}^n I(y^{(i)}=1)}\,
        (1-\psi)^{\sum\limits_{i=1}^n I(y^{(i)}=0)}d\psi\\[3pt]
  & = & U\textstyle
  \Big(f_0,\,f_1,\,\sum\limits_{i=1}^n\, I\big(y^{(i)}=0\big),\,
                   \sum\limits_{i=1}^n\, I\big(y^{(i)}=1\big)\Big)
\label{py-b}
\eeq
The function $U$ is defined as\vspace*{-5pt}
\beq
  U(f_0,f_1,n_0,n_1)
  & = & { \Gamma(f_0+f_1) \over \Gamma(f_0)\Gamma(f_1)} \,
           { \Gamma(f_0+n_0)\Gamma(f_1+n_1) \over \Gamma(f_0+f_1+n_0+n_1)}
  \ \ =\ \ { \prod\limits_{\ell=1}^{n_0} (f_0+\ell-1)\,
     \prod\limits_{\ell=1}^{n_1} (f_1+\ell-1)
     \over \prod\limits_{\ell=1}^{n_0+n_1} (f_0+f_1+\ell-1)}
  \ \ \ \ \ \ \
\eeq
The products above have the value one when the upper limits of $n_0$ or $n_1$ 
are zero.  The joint probability of $y\trn$ and the response, $y^*$, for
a test case is similar:
\beq
  P(y\trn,y^*) & = & U\textstyle
  \Big(f_0,\,f_1,\,\sum\limits_{i=1}^n\, I\big(y^{(i)}=0\big)\,
                  +I(y^*=0),\,
                  \sum\limits_{i=1}^n\, I\big(y^{(i)}=1\big)+I(y^*=1)\Big)
\label{pystar-b}
\eeq
Dividing $P(y\trn,y^*)$ by $P(y\trn)$ gives 
\beq
   P(y^*\ |\ y\trn) & = &
   \bern(y^*; \hat{\psi})
    \label{predp-y}
\eeq
Here, $\bern(y;\psi)\,=\,\psi^y\,(1-\psi)^{1-y}$ and
$\hat{\psi}\,=\,(f_1+N_1)\,/\,(f_0+f_1+n)$, with 
$N_y = \sum\limits_{\ell=1}^{n} I(y^{(\ell)}=y)$.\vspace*{-6pt}\linebreak
Note that $\hat{\psi}$ is just the posterior mean of $\psi$ based
on $y^{(1)},\ldots,y^{(n)}$.

Similarly, integrating over $\phi_{0,j}$ and $\phi_{1,j}$, we find 
that\vspace*{-3pt}
\beq
  P(x\trn_j\ |\ \theta_j,\,\alpha,\,y\trn) &\! =\! &
  \prod_{y=0}^1\,U(\alpha\theta_j,\,\alpha(1\!-\!\theta_j),\,
  O_{y,j},\,I_{y,j})
\label{pxj-b}\\[-13pt]\nonumber
\eeq
where $O_{y,j}\,=\,\sum\limits_{i=1}^n I(y^{(i)}=y,\,x_j^{(i)}=0)$ and
$I_{y,j}\,=\,\sum\limits_{i=1}^n I(y^{(i)}=y,\,x_j^{(i)}=1)$.

With $\psi$ and $\phib$ integrated out, we need deal only with
the remaining parameters, $\alpha$ and $\thetab$.  Note that
after eliminating $\psi$ and the $\phib$, the cases are no longer
independent (though they are exchangeable).  However, conditional 
on the responses, $y\trn$, and on $\alpha$, the values of
different features are still independent.  This is crucial to the
efficiency of the computations described below.


\subsection{\hspace*{-4pt}Predictions for test cases}\vspace*{-4pt}
\label{sub-pred-b}

We first describe how to predict the class for a test case when we are
either using all features, or using a subset of features without any
attempt to correct for selection bias.  We then consider how to make
predictions using our method of correcting for selection bias.

Suppose we wish to predict the response, $y^*$, in a test case for which we
know the retained features $\x_{1:k}^*=(\x^*_1,\cdots,\x^*_k)$ (having
renumbered features as necessary).  For this, we need the following
predictive probability:
\beq
  P(y^*\, |\, \x_{1:k}^*,\,\x_{1:k}\trn,\,y\trn) & = &
  { P(y^*\, |\, y\trn)\,P(\x_{1:k}^*\, |\, y^*,\,\x_{1:k}\trn,\,y\trn
    )_{\rule{0pt}{10pt}} \over
    \sum\limits_{y=0}^1
      P(y^*=y\, |\, y\trn)\,P(\x_{1:k}^*\, |\, y^*=y,\,\x_{1:k}\trn,\,y\trn)}
  \label{eq-unc-pp}
\eeq
Ie, we evaluate the numerator above for $y^*=0$ and
$y^*=1$, then divide by the sum to obtain the predictive probabilities.
The first factor in the numerator, $P(y^*\, |\, y\trn)$, is given
by equation~(\ref{predp-y}).  It is sufficient to obtain the
second factor up to a proportionality constant that doesn't depend on $y^*$,
as follows:
\beq
 P(\x_{1:k}^*\ |\ y^*,\,\x_{1:k}\trn,\,y\trn)
 &=& {P(\x_{1:k}^*,\,x_{1:k}\trn\ |\ y^*,\,y\trn)_{\rule{0pt}{9pt}} \over
     P(x_{1:k}\trn\ |\ y\trn)^{\rule{0pt}{7pt}}}
 \ \ \propto\ \ P(\x_{1:k}^*,\,x_{1:k}\trn\ |\ y^*,\,y\trn)
\eeq
This can be computed by integrating over $\alpha$, noting that conditional
on $\alpha$ the features are independent:\vspace*{-4pt}
\beq
  P(\x_{1:k}^*,\,x_{1:k}\trn\ |\ y^*,\,y\trn)
  & = & 
  \int P(\alpha)\,P(\x_{1:k}^*,\,x_{1:k}\trn\ |\ \alpha,\,y^*,\,y\trn)\,d\alpha
  \label{eq-alpha-int}\\[4pt]
  & = & 
  \int P(\alpha)\, \prod_{j=1}^k
   P(\x_j^*,\,x_j\trn\ |\ \alpha,\,y^*,\,y\trn)\,d\alpha
  \label{eq-alpha-int2}
\eeq
Each factor in the product above is found by using equation~(\ref{pxj-b})
and integrating over $\theta_j$:
\beq
  P(\x_j^*,\,x_j\trn\ |\ \alpha,\,y^*,\,y\trn) & = &
  \int_{0}^{1}\!
  P(\x_j^*\ |\ \theta_j,\,\alpha,\,\x_j\trn,\,y\trn,\, y^*)\,
  P(\x_j\trn\ |\ \theta_j,\,\alpha,\,y\trn)\,d\theta_j \ \ \ \ \ \ \\[2pt]
 &=&
 \int_{0}^{1}\!
  \bern(\x^*_j;\hat\phi_{y^*,j})\,
  \prod_{y=0}^1\,U(\alpha\theta_j,\,\alpha(1\!-\!\theta_j),
  \,I_{y,j},\,O_{y,j})\, d\theta_j\ \ \ \ \ \
  \label{eq-jfact}
\eeq
where $\hat\phi_{y^*,j}= (\alpha\theta_j+I_{y^*,j})\ /\
(\alpha+N_{y^*})$, the posterior mean of $\phi_{y^*,j}$ given $\alpha$
and $\theta_j$. 

When using $k$ features selected from a larger number, $p$, the predictions
above, which are conditional on only $x\trn_{1:k}$ and $y\trn$, are not correct
--- we should also condition on the event,
$\mathcal{S}$, that $|\COR (y\trn,x\trn_j)|\,\le\,\gamma$ for $j=k+1,\ldots,p$.
We need to modify the predictive probability of equation~(\ref{eq-unc-pp})
by replacing $P(\x_{1:k}^*\ |\ y^*,\,\x_{1:k}\trn,\,y\trn)$ with
$P(\x_{1:k}^*\ |\ y^*,\,\x_{1:k}\trn,\,y\trn,\,\mathcal{S})$, which is
proportional to $P(\x_{1:k}^*,\,\x_{1:k}\trn,\,\mathcal{S}\ |\ y^*,\,y\trn)$.
Analogously to equations~(\ref{eq-alpha-int}) and~(\ref{eq-alpha-int2}), 
we obtain
\beq
  P(\x_{1:k}^*,\,x_{1:k}\trn,\,\mathcal{S}\ |\ y^*,\,y\trn)
  & = & 
  \int P(\alpha)\,P(\x_{1:k}^*,\,x_{1:k}\trn,\,\mathcal{S}
    \ |\ \alpha,\,y^*,\,y\trn)\,d\alpha
  \label{eq-alpha-int-mod}\\[4pt]
  & = & 
  \int P(\alpha)\, P(\mathcal{S}\ |\ \alpha,\,y\trn) \prod_{j=1}^k
   P(\x_j^*,\,x_j\trn\ |\ \alpha,\,y^*,\,y\trn)\,d\alpha\ \ \ \
  \label{eq-alpha-int2-mod}
\eeq
The factors for the $k$ retained features are computed as before, using
equation~(\ref{eq-jfact}).  The additional correction factor that is needed
(presented earlier as equation~(\ref{eq-cor-lik})) is
\beq
  P(\mathcal{S}\ |\ \alpha,\,y\trn) & = &
    \prod_{j=k+1}^p P(|\COR (y\trn,x\trn_j)|\,\le\,\gamma\ |\ \alpha,\,y\trn)
  \\[3pt]
  & = &
  \Big[\, 
      P(|\COR (y\trn,x\trn_t)|\,\le\,\gamma\ |\ \alpha,\,y\trn)
    \,\Big]^{p-k}
  \label{eq-corfact}
\eeq
where $t$ is any of the omitted features, all of which have the same
probability of having a small correlation with $y$.  
We discuss how to compute this adjustment factor in the next section.

To see intuitively why this adjustment factor will correct for
selection bias, recall that as discussed in
Section~(\ref{sub-model-b}), when $\alpha$ is small, features will be
more likely to have a strong relationship with the response. If the
likelihood of $\alpha$ is based only on the selected features, which
have shown high correlations with the response in the training
dataset, it will favor values of $\alpha$ that are inappropriately
small.  Multiplying by the adjustment factor, which favors larger
values for $\alpha$, undoes this bias.

We compute the integrals over $\alpha$ in equations~(\ref{eq-alpha-int2})
and~(\ref{eq-alpha-int2-mod}) by numerical quadrature.  We use the midpoint
rule, applied to $u=F^{-1}(\alpha)$, where $F^{-1}$ is the inverse cumulative
distribution function for the Inverse-Gamma$(a,b)$ prior for $\alpha$. The
prior for $u$ is uniform over $(0,1)$, and so needn't be explicitly included
in the integrand.  With $K$ points for the midpoint rule, the effect is that
we average the value of the integrand, without the prior factor, for values of
$\alpha$ that are the $0.5/K, 1.5/K, \ldots, 1-0.5/K$ quantiles of its
Inverse-Gamma prior. For each $\alpha$, we use Simpson's Rule to compute the
one-dimensional integrals over $\theta_j$ in equation~(\ref{eq-jfact}).  

\subsection{\hspace*{-4pt}Computation of the adjustment factor}\vspace*{-4pt}
\label{sec-adj-b}

Our remaining task is to compute the adjustment factor of
equation~(\ref{eq-corfact}), which depends on the
probability that a feature will have correlation less than $\gamma$ in
absolute value.  Computing this seems difficult ---
we need to sum the probabilities of $\x_t\trn$ given $y\trn$,
$\alpha$ and $\theta_t$ over all configurations of $\x_t\trn$ 
for which $|\COR (y\trn,x\trn_t)|\,\le\,\gamma$ ---
but the computation can be 
simplified by noticing that $\COR(x\trn_t,y\trn)$ can be written in terms of
\mbox{$I_0\,=\,\sum_{i=1}^n I(y^{(i)}=0,\,x^{(i)}_t=1)$} and
$I_1\,=\,\sum_{i=1}^n I(y^{(i)}=1,\,x^{(i)}_t=1)$, as follows:\vspace*{-7pt}
\beq
  \COR(x_t\trn,\,y\trn)
  & = &
  {\displaystyle\sum_{i=1}^n\, \big(y^{(i)}-\bar y\big)\,x^{(i)}_t
  \over
   \sqrt{\sum\limits_{i=1}^n \big(y^{(i)}-\bar y\big)^2}\
   \sqrt{\sum\limits_{i=1}^n \big(x^{(i)}_t-\bar x_t\big)^2} } \\[6pt]
  & = &
  { (0-\ybar)\,I_0\ +\ (1-\ybar)\,I_1 \over
     \sqrt{n \ybar (1\!-\!\ybar)}\,
     \sqrt{I_0+I_1-(I_0+I_1)^2/n} }
  \label{cor-express}
\eeq
We write the above as $\Cor(I_0,I_1,\ybar)$, taking $n$ as known.
This function is defined for $0 \le I_0 \le n(1\!-\!\ybar)$ and
$0 \le I_1 \le n\ybar$.

\begin{figure}

\hspace*{0.2in}\psfig{file=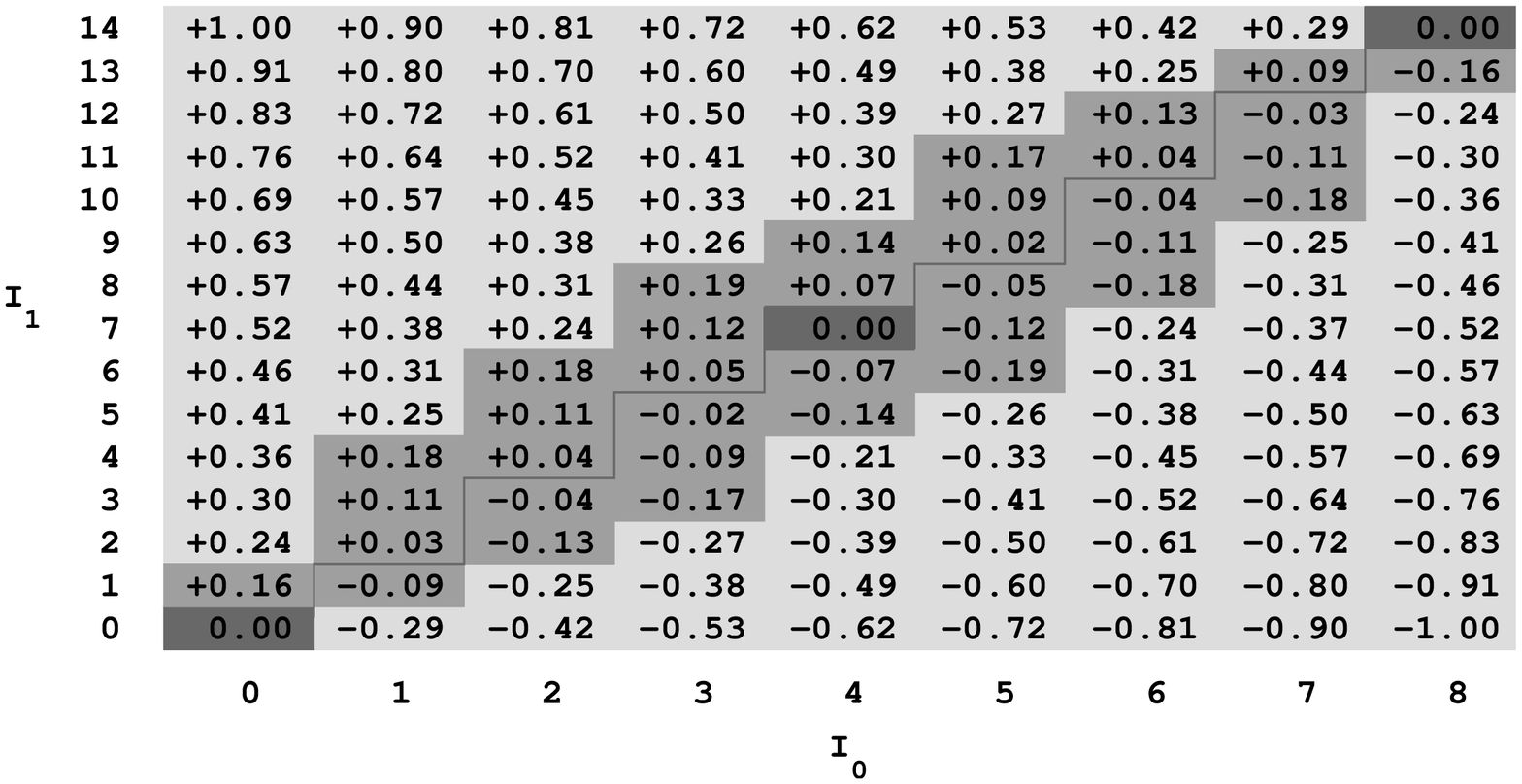,width=6in}

\caption[]{The $\Cor$ function for a dataset with $n=22$ and $\ybar=14/22$.
The values of $\Cor(I_0,I_1,\ybar)$ are shown for the valid range of
$I_0$ and $I_1$.  Using $\gamma=0.2$, the values of ($I_0$,$I_1$) in
$L_0$ are shown in dark
grey, those in $L_-$ or $L_+$ in medium grey, and those in $H_-$ or $H_+$
in light grey.}\label{fig-sets}

\end{figure}

Fixing $n$, $\ybar$, and $\gamma$, we can define the following sets
of values for $I_0$ and $I_1$ (for some feature $x_t$) in terms of the
resulting correlation with $y$:
\beq
   L_0 & = & \{\,(I_0,I_1)\ :\ \Cor(I_0,I_1,\ybar) = 0 \,\} \\[4pt]
   L_+ & = & \{\,(I_0,I_1)\ :\ 0 < \Cor(I_0,I_1,\ybar) \le \gamma \,\} \\[4pt]
   L_- & = & \{\,(I_0,I_1)\ :\ -\gamma \le \Cor(I_0,I_1,\ybar) < 0 \,\} \\[4pt]
   H_+ & = & \{\,(I_0,I_1)\ :\ \gamma < \Cor(I_0,I_1,\ybar) \,\} \\[4pt]
   H_- & = & \{\,(I_0,I_1)\ :\ \Cor(I_0,I_1,\ybar) < -\gamma \,\}
\eeq
A feature will be discarded if $(I_0,I_1) \,\in\, L_- \cup L_0 \cup L_+$
and retained if $(I_0,I_1) \,\in\, H_- \cup H_+$.  These sets are
illustrated in Figure~\ref{fig-sets}.

We can write the probability needed in equation~(\ref{eq-corfact}) using
either $L_-$, $L_0$, and $L_+$ or $H_-$ and $H_+$.  We will take the latter 
approach here, as follows:
\beq
\ \ \ \
P(\, |\COR(x\trn_t,y\trn)| \leq \gamma \ |\ \alpha,\,y\trn)
 & = &
   1\ -\ P(\,(I_0,I_1)\,\in\, H_- \cup H+\ |\ \alpha,\,y\trn) \\[6pt]
 & = & 1\ -\ \!\!\!\!\!\!\sum_{
   \begin{array}{cc} \\[-19pt] \scriptstyle (I_0,I_1)\, \in \\[-3pt]
                     \scriptstyle H_- \cup H_+ \\[-8pt]
   \end{array}}\!\!\!\!\!\!
   P(I_0,\,I_1\ |\ \alpha,\,y\trn)
 \label{adj-fact-H1}
\eeq

We can now exploit symmetries of the prior and of the $\Cor$ function
to speed up computation.  First, note that
$\Cor(I_0,I_1,\ybar)\, =\, -\Cor(n(1\!-\!\ybar)-I_0,n\ybar-I_1,\ybar)$,
as can be derived from equation~(\ref{cor-express}), or by simply 
noting that swapping the feature values (0 and 1) should change only the sign
of the correlation.  The one-to-one mapping $(I_0,I_1) \rightarrow
(n(1\!-\!\ybar)-I_0,n\ybar-I_1)$, which maps $H_-$ and $H_+$ and vice
versa (similarly for $L_-$ and $L_+$), therefore leaves $\Cor$ unchanged.
The priors for
$\theta$ and $\phi$ (see~(\ref{prior-theta-b}) and~(\ref{prior-phi-b}))
are symmetrical with respect to the class labels 0 and 1, so the prior
probability of
($I_0,\,I_1)$ is the same as that of $(n(1\!-\!\ybar)-I_0,\,n\ybar-I_1)$.
We can therefore rewrite equation~(\ref{adj-fact-H1}) as
\beq
P(\, |\COR(x\trn_t,y\trn)| \leq \gamma \ |\ \alpha,\,y\trn) & = & 1\ -\
      2\!\!\!\!\!\!\sum_{(I_0,I_1)\, \in\, H_+}\!\!\!\!\!\!\!
      P(I_0,\,I_1\ |\ \alpha,\,y\trn)
 \label{adj-fact-H2}
\eeq

At this point we write the probabilities for $I_0$ and $I_1$ in
terms of an integral over $\theta_t$, and then swap the order of summation and
integration, obtaining
\beq
  \sum_{(I_0,I_1)\, \in\, H_+}\!\!\!\!\!\!\!
      P(I_0,\,I_1\ |\ \alpha,\,y\trn) & = &
  \int_{0}^{1}\!\!
  \sum_{(I_0,I_1)\, \in\, H_+}\!\!\!\!\!
    P(I_0,\,I_1\ |\ \alpha,\,\theta_t,\,y\trn)\ d\theta_t
\label{eq-adjint}
\eeq
The integral over $\theta_t$ can be approximated using
some one-dimensional numerical quadrature method (we use Simpson's Rule),
provided we can evaluate the integrand.

The sum over $H_+$ can easily be delineated because
$\Cor(I_0,I_1,\ybar)$ is a monotonically decreasing function of $I_0$,
and a monotonically increasing function of $I_1$, as may be confirmed
by differentiating with respect to $I_0$ and $I_1$.  Let $b_0$ be the
smallest value of $I_1$ for which $\Cor(0,I_1,\ybar) > \gamma$.  Taking
the ceiling of the solution of $\Cor(0,I_1,\ybar)=\gamma$, we find
that $b_0\, =\,\lceil 1/(1/n+(1-\bar{y})/(n\bar{y}\gamma^2))\rceil$.
For $b_0 \le I_1 \le n\ybar$, let $r_{I_1}$ be the largest value of
$I_0$ for which $\Cor(I_0,I_1,\ybar) > \gamma$.  We can write
\beq
  \sum_{(I_0,I_1)\, \in\, H_+}\!\!\!\!\!
    P(I_0,\,I_1\ |\ \alpha,\,\theta_t,\,y\trn) & = &
  \sum_{I_1=b_0}^{n\ybar}\, \sum_{I_0=0}^{r_{I_1}}\,
    P(I_0,\,I_1\ |\ \alpha,\,\theta_t,\,y\trn)
\eeq

Given $\alpha$ and $\theta_t$, $I_0$ and $I_1$ are independent, so
we can reduce the computation needed by rewriting the above expression
as follows:
\beq
  \sum_{(I_0,I_1)\, \in\, H_+}\!\!\!\!\!
    P(I_0,\,I_1\ |\ \alpha,\,\theta_t,\,y\trn) & = &
  \sum_{I_1=b_0}^{n\ybar}\,P(I_1\ |\ \alpha,\,\theta_t,\,y\trn)\,
  \sum_{I_0=0}^{r_{I_1}}\, P(I_0\ |\ \alpha,\,\theta_t,\,y\trn)
\label{eq-sum-H}
\eeq
Note that the inner sum can be
updated from one value of $I_1$ to the next by just adding any additional
terms needed.
This calculation therefore requires $1\!+\!n\ybar\!-\!b_0 \le n$ 
evaluations of $P(I_1\ |\
\alpha,\,\theta_t,\,y\trn)$ and $1\!+\!r_{n\ybar} \le n$ evaluations of
$P(I_0\ |\ \alpha,\,\theta_t,\,y\trn)$.  

To compute $P(I_1\ |\ \alpha,\,\theta_t,\,y\trn)$, we multiply
the probability of any particular value for $x\trn_t$ in which there are $I_1$
cases with $y=1$ and $x_t=1$ by the number of ways this can occur.
The probabilities are found 
by integrating
over $\phi_{0,t}$ and $\phi_{1,t}$, as described in Section~\ref{sub-int-b}.
The result is\vspace*{-6pt}
\beq
  P(I_1\ |\ \alpha,\,\theta_t,\,y\trn) & = &
    \mychoose{n\ybar}{I_1} U(\alpha\theta_t,\,\alpha(1\!-\!\theta_t),\,
                            I_1,\,n\ybar-I_1)
\eeq
Similarly,
\beq
  P(I_0\ |\ \alpha,\,\theta_t,\,y\trn) & = &
    \mychoose{n(1\!-\!\ybar)}{I_0} U(\alpha\theta_t,\,\alpha(1\!-\!\theta_t),\,
                            I_0,\,n(1\!-\!\ybar)-I_0)
\eeq
One can easily derive simple expressions for
$P(I_1\ |\ \alpha,\,\theta_t,\,y\trn)$ and
$P(I_0\ |\ \alpha,\,\theta_t,\,y\trn)$ in
terms of $P(I_1-1\ |\ \alpha,\,\theta_t,\,y\trn)$ and
$P(I_0-1\ |\ \alpha,\,\theta_t,\,y\trn)$, which avoid the need to compute
gamma functions or large products for each value of $I_0$ or $I_1$
when these values are used sequentially, as in equation~(\ref{eq-sum-H}).

\section{\hspace*{-7pt}A Simulation Experiment
   }\vspace*{-8pt}\label{sec-sim-bnaive}

In this section, we use a dataset generated from the naive Bayes model
defined in Section~\ref{sub-model-b} to demonstrate the lack of
calibration that results when only a subset of features is used,
without correcting for selection bias.  We show that our
bias-correction method eliminates this lack of calibration.  We will
also see that for the naive Bayes model only a small amount of extra
computation time is needed to obtain the adjustment factor needed
by our method.

\begin{figure}[p]

\centerline{\includegraphics{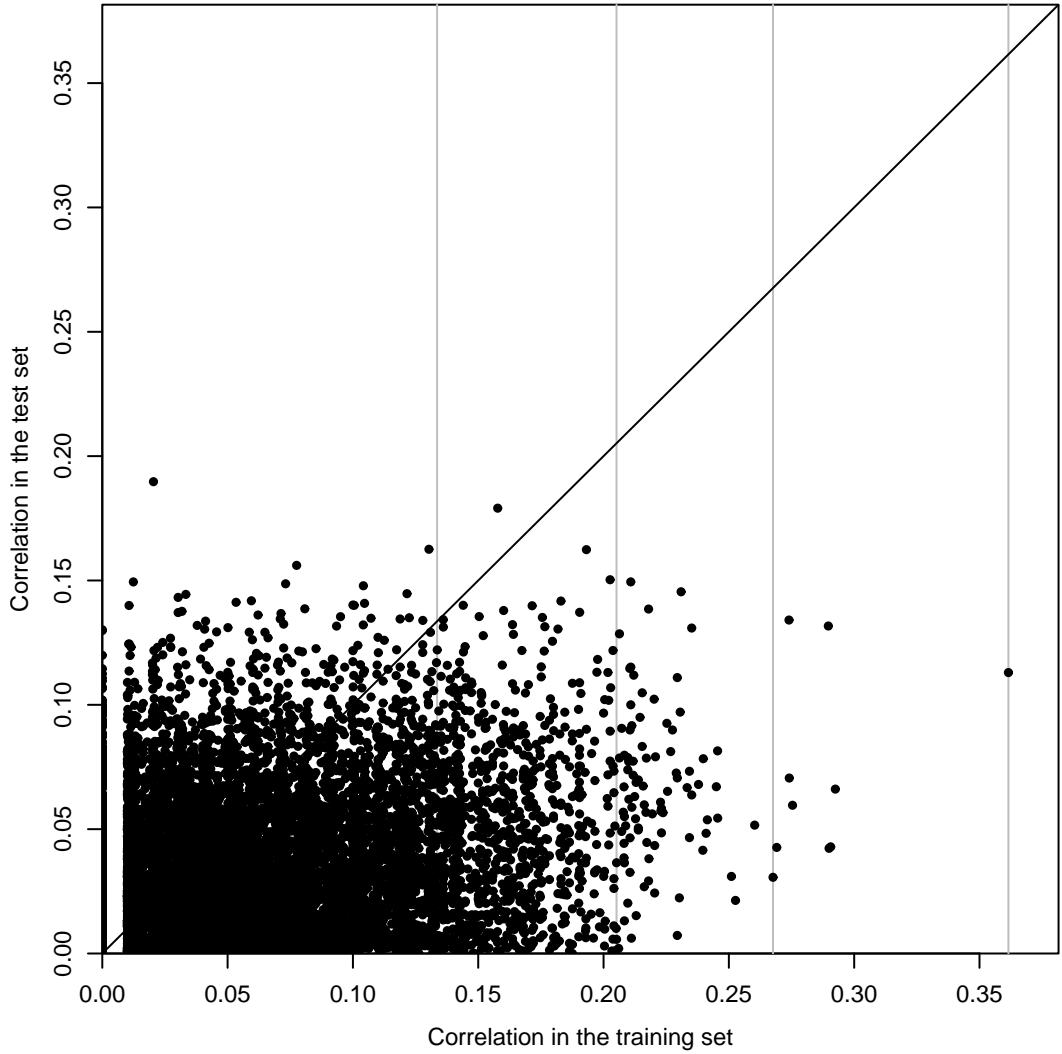}}

\caption[]{The absolute value of the sample correlation of each
feature with the binary response, in the training set, and in the test
set.  Each dot represents one of the 10000 binary features. The
training set correlations of the 1st, 10th, 100th, and 1000th most
correlated features are marked by vertical
lines.}\label{fig-trn-tst-cor}

\end{figure}

Fixing $\alpha=300$, and $p=10000$, we used equations
(\ref{sample-x-b}), (\ref{prior-theta-b}) and (\ref{prior-phi-b}) to
generate a set of 200 training cases and a set of 2000 test cases,
both having equal numbers of cases with $y=0$ and $y=1$.  We then
selected four subsets of features, containing 1, 10, 100, and 1000
features, based on the absolute values of the sample correlations of
the features with $y$.  The smallest correlation (in absolute value)
of a selected feature with the class was 0.36, 0.27, 0.21, and 0.13
for these four subsets.  These are the values of $\gamma$ used by the
bias correction method when computing the adjustment factor of
equation~(\ref{eq-corfact}).  Figure~\ref{fig-trn-tst-cor} shows the
absolute value of the sample correlation in the training set of all
10000 features, plotted against the sample correlation in the test
set.  As can be seen, the high sample correlation of many selected
features in the training set is partly or wholely a matter of chance,
with the sample correlation in the test set (which is close to the
real correlation) often being much less.  The role of chance is
further illustrated by the fact that the feature with highest sample
correlation in the test set is not even in the top 1000 by sample
correlation in the training set.

For each number of selected features, we fit this data using the naive
Bayes model with the prior for $\psi$ (equation~(\ref{prior-psi-b}))
having $f_0=f_1=1$ and the prior for $\alpha$
(equation~(\ref{prior-alpha-b})) having shape parameter $a=0.5$ and
rate parameter $b=5$.  We then made predictions for the test cases
using the methods described in Section~\ref{sub-pred-b}.  The
``uncorrected'' method, based on equation~(\ref{eq-unc-pp}), makes no
attempt to correct for the selection bias, whereas the ``corrected''
method, with the modification of equation~(\ref{eq-alpha-int2-mod}),
produces predictions that account for the procedure used to select the
subset of features.  We also made predictions using all 10000
features, for which bias correction is unnecessary.

\begin{figure}[p]

\includegraphics{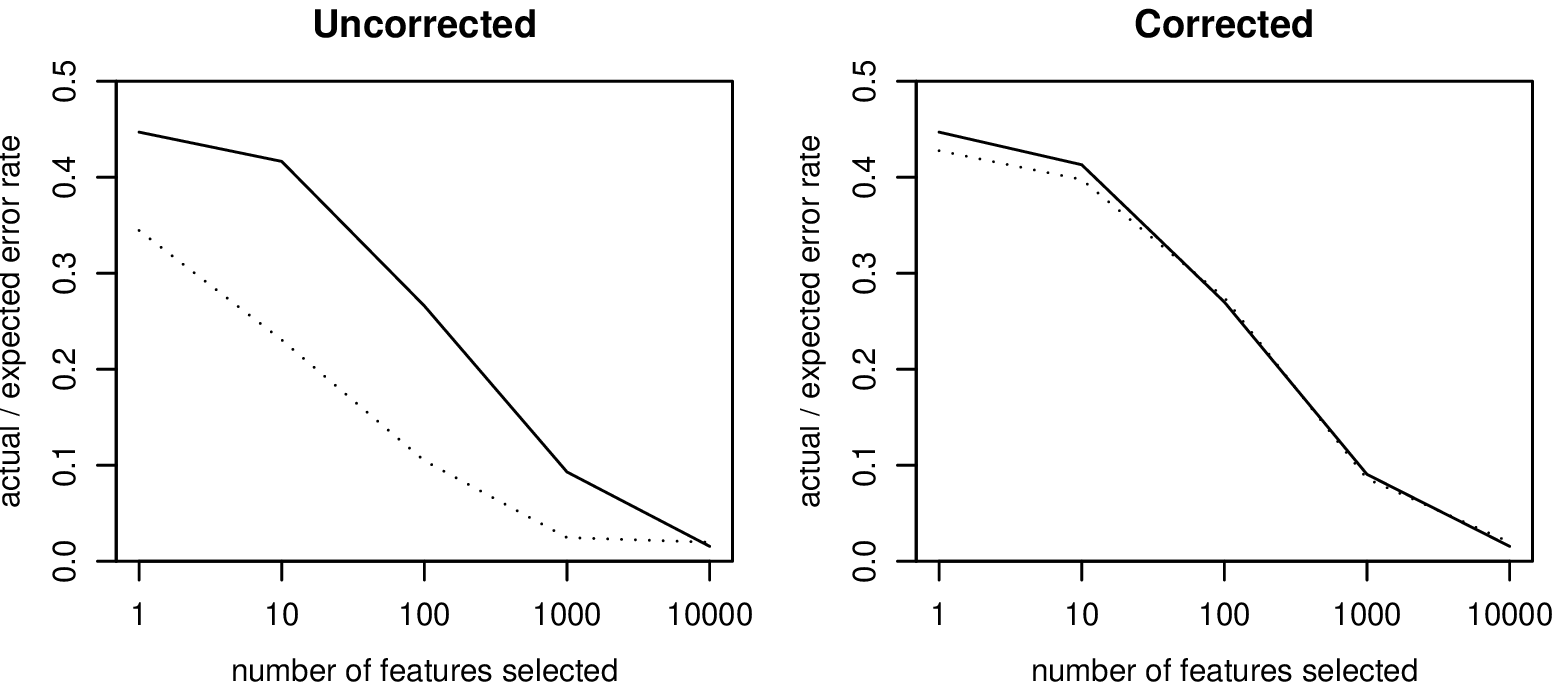}

\vspace*{-6pt}

\caption[]{Actual and expected error rates with varying numbers of 
features selected, with and without correction for selection bias.
The solid line is the actual error rate on test cases. The dotted
line is the error rate that would be expected based on the 
predictive probabilities.}\label{fig-sim-perf1}

\end{figure}

We compared the predictive performance of the corrected method with
the uncorrected method in several ways.  First, we looked at the error
rate when classifying test cases by thresholding the predictive
probabilities at $1/2$.  As can be seen in Figure~\ref{fig-sim-perf1},
there is little difference in the error rates with and without
correction for bias.  However, the methods differ drastically in terms
of the \textit{expected} error rate --- the error rate we would expect
based on the predictive probabilities for the test cases, equal to
$(1/N) \sum_i\, \hat p^{(i)}\, I(\hat p^{(i)}<0.5) \ +\ (1\!-\!\hat
p^{(i)})\, I(\hat p^{(i)}\ge0.5)$, where $\hat p^{(i)}$ is the
predictive probability of class 1 for test case $i$.  The predictive
probabilities produced by the uncorrected method would lead us to
believe that we would have a much lower error rate than the actual
performance.  In contrast, the expected error rates based on the
predictive probabilities produced using bias correction closely
match the actual error rates.

\begin{figure}[p]

\vspace*{-12pt}

\includegraphics{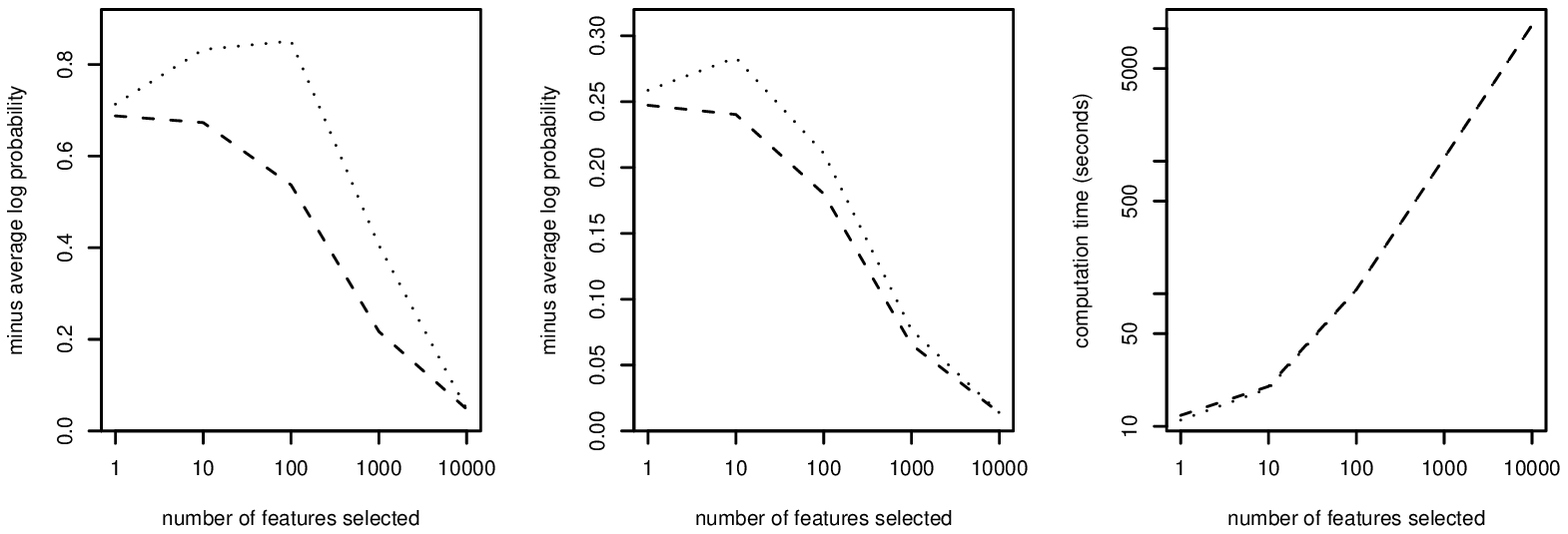}

\vspace*{-6pt}

\caption[]{Performance in terms of average minus log probability and
average squared error, and computation time, with varying
numbers of features selected, with and without correction for selection bias.
The left plot shows minus the average log probability of the correct class
for test cases, with 1, 10, 100, 1000, and all 10000 features selected.
The dashed line is with bias correction, the dotted line without.  The
middle plot is similar, but shows average squared error on test cases. 
The right plot shows the computation time needed for the methods (the
two lines almost coincide).  Note
that when all 10000 features are used, there is no difference between the
corrected and uncorrected methods.}\label{fig-sim-perf2}

\end{figure}

Two additional measures of predictive performance are shown in
Figure~\ref{fig-sim-perf2}.  One measure of performance is minus the
average log probability of the correct class in the $N$ test cases,
which is $-(1/N)\,\sum_{i=1}^N\,
[y^{(i)}\log(\hat{p}^{(i)})\,+\,(1\!-\!y^{(i)})\log(1\!-\!\hat{p}^{(i)})]$.
This measure heavily penalizes test cases where the actual class has a
predictive probability near zero.  Another measure, less sensitive to
such drastic errors, is the average squared error between the actual
class (0 or 1) and the probability of class 1, given by $(1/N)
\sum_{i=1}^N (y^{(i)}-\hat p ^{(i)})^2$.  The corrected method
outperforms the uncorrected method by both these measures, with the
difference being greater for minus average log probability.
Interestingly, performance of the uncorrected method actually gets
worse when going from 1 feature to 10 features.  This may be because
the single feature with highest sample correlation with the response
does have a strong relationship with the response (as may be likely in
general), whereas some other of the top 10 features by sample
correlation have little or no real relationship.

\begin{table}
\small \input{tab}

\vspace*{-2.02in}
\hspace*{2.4in}\begin{minipage}{4.0in}\caption[]{
Comparison of calibration for predictions found with and without
correction for selection bias, on data simulated from the binary naive 
Bayes model.
Results are shown with four subsets of features and with the
complete data (for which no correction is necessary). 
The test cases were divided into 10 categories by the first decimal of the
predictive probablity of class 1.  The table shows the number 
of test cases in each category for each method (``\#''), the
average predictive probability of class 1 for cases in that category 
(``Pred''), and the actual fraction of these cases that were in class 1
(``Actual'').}\label{tab-bnaive}
\end{minipage}

\end{table}

We also looked in more detail at how well calibrated the predictive
probabilities were.  Table~\ref{tab-bnaive} shows the average
predictive probability for class 1 and the actual fraction of cases in
class 1 for test cases grouped according to the first
decimal of their predictive probabilities, for both the uncorrected
and the corrected methods.  Results are shown using subsets of
1, 10, 100, and 1000 features, and using all features.
We see that the uncorrected method produces overconfident
predictive probabilities, either too close to zero or
too close to one.  The corrected method avoids such bias (the values
for ``Pred'' and ``Actual'' are much closer), showing that it is well
calibrated.  

\begin{figure}[t]
\vspace*{-6pt}

\hspace*{12pt}\includegraphics[width=6.4in]{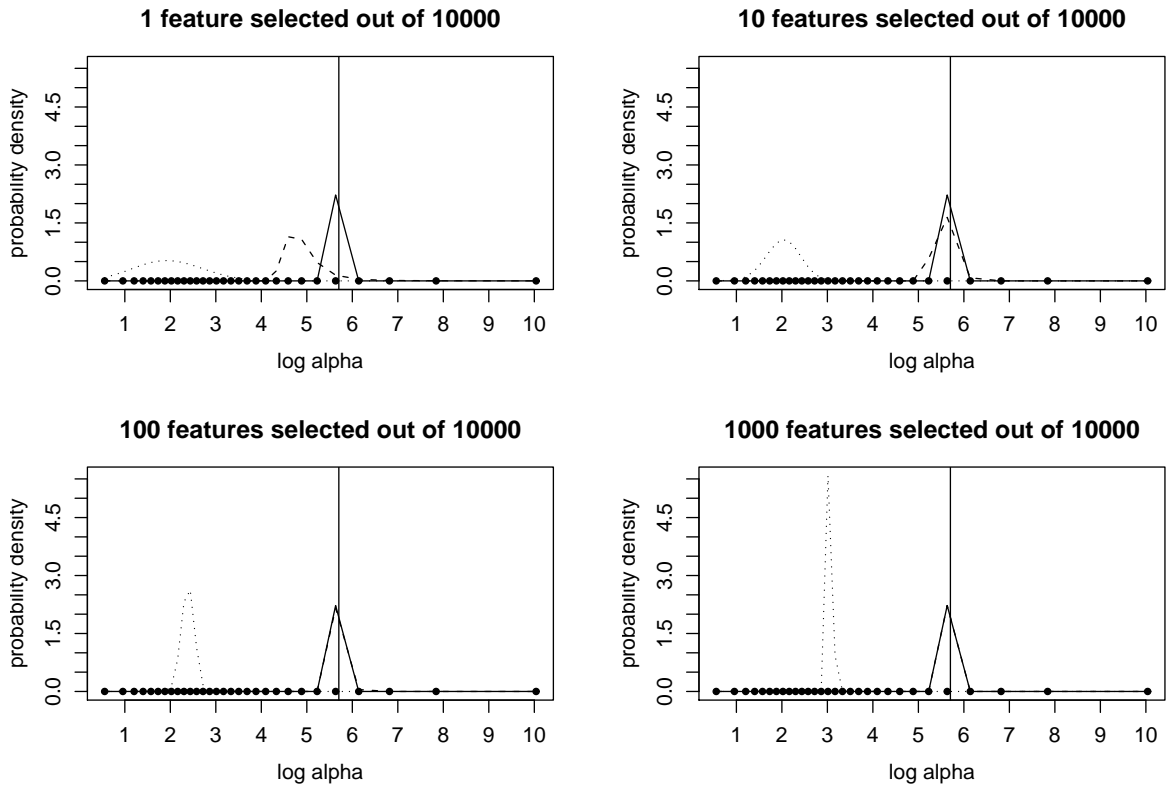}

\vspace*{-10pt}

\caption{Posterior distributions of $\log(\alpha)$ for the
simulated data, with different numbers of features selected.
The true value of $\log(\alpha)$ is 5.7, shown by the vertical line.
The solid line is the posterior density using all features.  For each number
of selected features, the dashed line is the posterior density
including the factor that corrects for selection bias; the dotted line is
the posterior density without bias correction.  
The dashed and solid lines overlap in the bottom two graphs.
The dots mark the values of $\log(\alpha)$ used to approximate the density, at
the $0.5/K,1.5/K,\ldots,(K\!-\!0.5)/K$ quantiles of the prior distribution
(where $K=30$).
The probabilities of $x\trn$ at each of these values for 
$\alpha$ were computed, rescaled to sum to $K$, and finally 
multiplied by the Jacobian, $\alpha P(\alpha)$,
to obtain the approximation to the posterior density of 
$\log(\alpha)$.
}
\label{fig-post-alpha}
\end{figure}

The biased predictions of the uncorrected method result from an
incorrect posterior distribution for $\alpha$, as illustrated in
Figure~\ref{fig-post-alpha}.  Without bias correction, the posterior
based on only the selected features incorrectly favours values of
$\alpha$ smaller than the true value of 300.  Multiplying by the
adjustment factor corrects this bias in the posterior distribution.

Our software (available from
\texttt{http://www.utstat.utoronto.ca/$\sim$longhai}) is written in the
R language, with some functions for intensive computations such as
numerical integration and computation of the adjustment factor written
in C for speed.  We approximated the integral with respect to $\alpha$
using the midpoint rule with $K=30$ values for $F^{-1}(\alpha)$, as
discussed at the end of Section~\ref{sub-pred-b}.  The integrals with
respect to $\theta$ in equations~(\ref{eq-jfact})
and~(\ref{eq-adjint}) were approximated using Simpson's Rule,
evaluating $\theta$ at 21 points.

Computation times for each method (on a 1.2 GHz UltraSPARC III
processor) are shown on the right in Figure~\ref{fig-sim-perf2}.  The
corrected method is almost as fast as the uncorrected method, since
the time to compute the adjustment factor is negligible compared to
the time spent computing the integrals over $\theta_j$ for the
selected features.  Accordingly, considerable time can be saved by
selecting a subset of features, rather than using all of them, without
introducing an optimistic bias, though some accuracy in predictions
may of course be lost when we discard the information contained in the
unselected features.

\section{\hspace*{-7pt}A test using gene expression data}\vspace*{-8pt}
\label{sec-gene}

We also tested our method using a publicly available dataset on gene
expression in normal and cancerous human colon tissue.  This dataset
contains the expression levels of 6500 genes in 40 cancerous and 22 normal
colon tissues, measured using the Affymetrix technology.  The dataset
is available at
\texttt{http://geneexpression.cinj.org/$\sim$notterman/affyindex.html}.
We used only the 2000 genes with highest minimal intensity, as
selected by Alon, Barkai, Notterman, Gish, Mack, and Levine (1999).
In order to apply the binary naive Bayes model to the data, we
transformed the real-value data into binary data by thresholding at
the median, separately for each feature.

We divided these 2000 genes randomly into 10 equal groups, producing
10 smaller datasets, each with 200 binary features, as well as the
binary class (normal/cancerous).  We applied the corrected and
uncorrected methods separately to each of these 10 datasets, allowing
some assessment of variability when comparing performance.  For each
of these 10 datasets, we used leave-one-out cross validation to obtain
predictive probabilities for the class in the 62 cases.  In this
cross-validation procedure, we left out each of the 62 cases in turn,
selected the five features with the largest sample correlation with
the class (in absolute value), and found the predictive probability
for the left-out case using the binary naive Bayes model, with and
without bias correction.  The absolute value of the correlation of the
last selected feature with the class was always around 0.5.  We used
the same prior distribution, and the same computational methods, as
for the demonstration in Section~\ref{sec-sim-bnaive}.

\begin{figure}[p]

\hspace*{-0.1in}\includegraphics{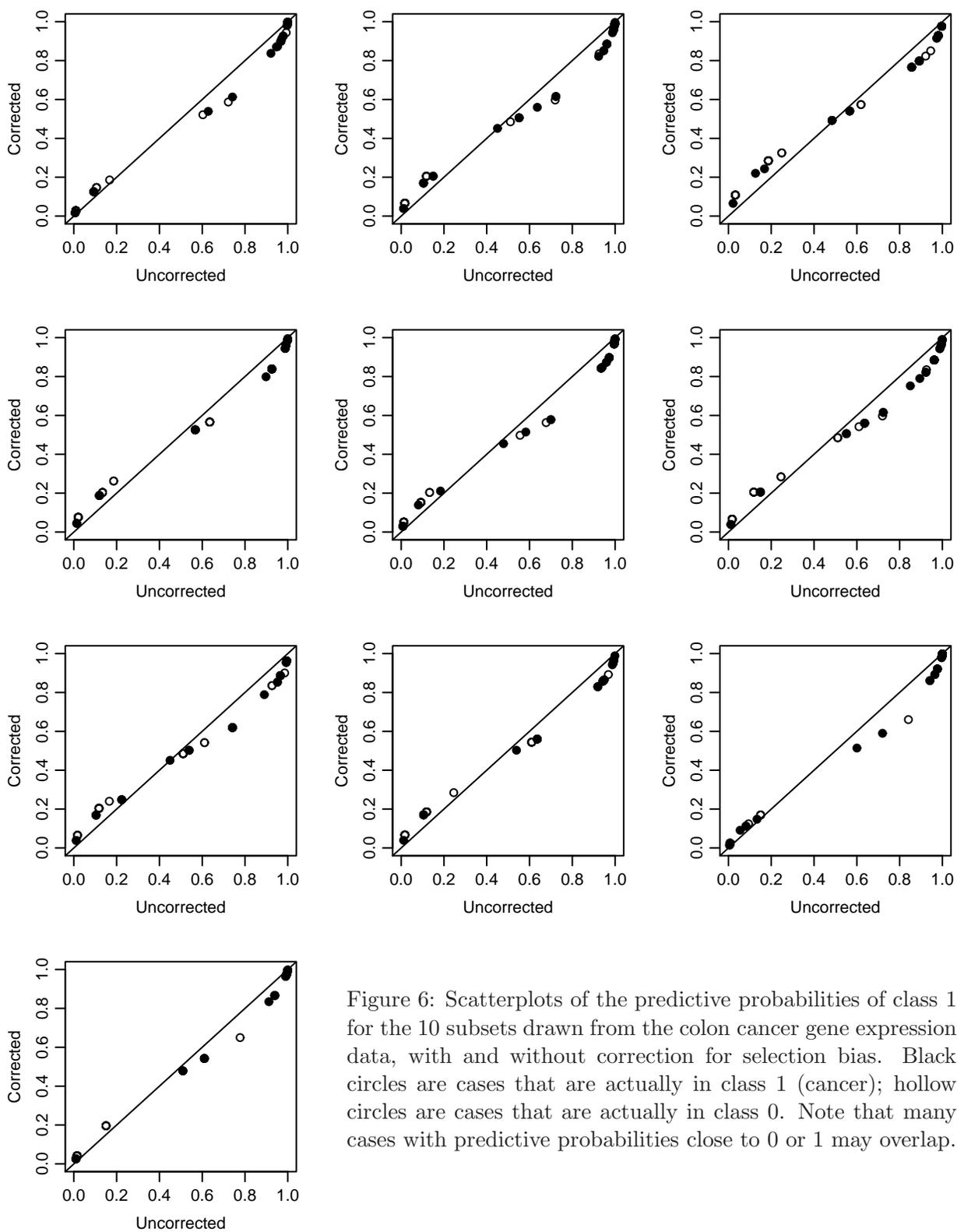}

\vspace*{-2in}

\hspace*{2.3in}\begin{minipage}{4.1in}\caption{Scatterplots of the
predictive probabilities of class 1 for the 10 subsets drawn from
the colon cancer gene expression data, with and without correction for
selection bias.  Black circles are cases that are actually in class 1
(cancer); hollow circles are cases that are actually in class 0. Note
that many cases with predictive probabilities close to 0 or 1 may 
overlap.}\label{fig-colon-pred-prob}\end{minipage}

\vspace*{0.4in}

\end{figure}

Figure~\ref{fig-colon-pred-prob} plots the predictive probabilities of
class 1 for all cases, with each of the 10 subsets of features.  The
tendency of the uncorrected method to produce more extreme
probabilities (closer to 0 and 1) is clear.  However, when the
predictive probability is close to 0.5, there is little difference
between the corrected and uncorrected methods.  Accordingly, the two
methods usually classify cases the same way, if classification is done
by thresholding the predictive probability at 0.5, and have very
similar error rates. (The overall average error rate is 0.194 for the
uncorrected method and 0.182 for the corrected method.)  Note,
however, that correcting for bias would have a substantial effect if
cases were classified by thresholding the predictive probability at
some value other than 0.5, as would be appropriate if the consequences
of an error are different for the two classes.

%
%

\begin{figure}[p] 

\begin{center}

\vspace*{-0.2in}

\includegraphics{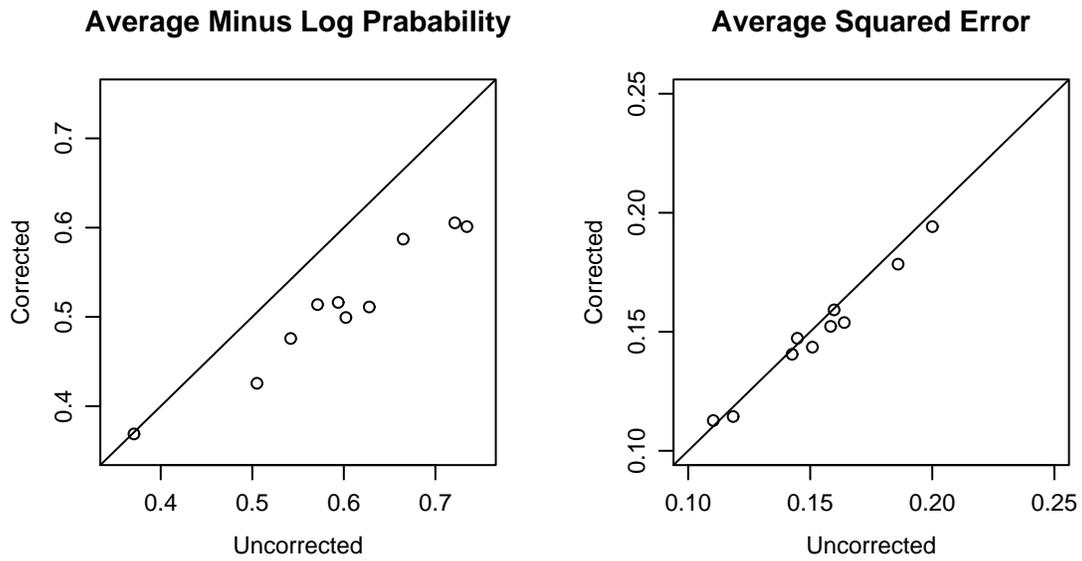}

\vspace*{-0.25in}

\end{center} 

\caption{Scatterplots of the average minus log probability of the
correct class and of the average squared error (assessed by cross
validation) when using the 10 subsets of features for the colon cancer
gene expression data, with and without correcting for selection bias.}
\label{fig-colon}

\end{figure}

\begin{figure}[p] 

\begin{center}

\vspace*{-0.2in}

\includegraphics{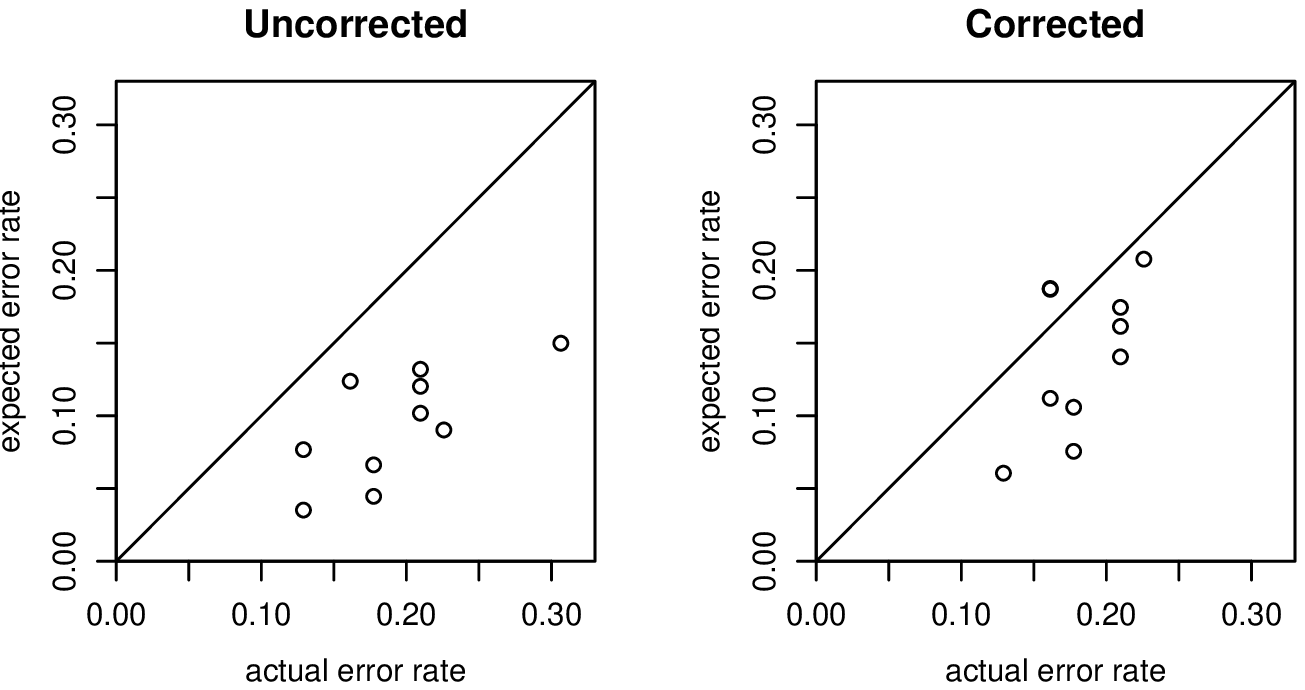}

\vspace*{-0.25in}

\end{center} 

\caption[]{Actual versus expected error rates on the colon cancer
datasets, with and without bias correction.  Points are shown for each
of the 10 subsets of features used for testing.  In the plot on the
right, two points above the diagonal are almost superimposed.
} \label{fig-colon2}

\end{figure}

Figure~\ref{fig-colon} compares the two methods in terms of average
minus log probability of the correct class and in terms of average
squared error.  From these plots it is clear that bias correction
improves the predictive probabilities.  In terms of average minus log
probability, the corrected method is better for all 10 datasets, and
in terms of average squared error, the corrected method is better for
8 out of 10 datasets.  (A paired $t$ test with these two measures
produced $p$-values of $0.00007$ and $0.019$ respectively.)

Finally, Figure~\ref{fig-colon2} shows that our bias correction method
reduces optimistic bias in the predictions.  For each of the 10
datasets, this plot shows the actual error rate (in the leave-one-out
cross-validation assessment) and the error rate expected from the
predictive probabilities.  For all ten datasets, the expected error
rate with the uncorrected method is substantially less than the actual
error rate.  This optimistic bias is reduced in the corrected method,
though it is not eliminated entirely.  The remaining bias presumably
results from the failure in this dataset of the naive Bayes assumption
that features are independent within a class.

\section{\hspace*{-7pt}Conclusions and Future Work}
\label{sec-conc}\vspace*{-8pt}

We have proposed a Bayesian method for making well-calibrated
predictions for a response variable when using a subset of features
selected from a larger number based on some measure of dependency
between the feature and the response.  Our method results from
applying the basic principle that predictive probabilities should be
conditional on all available information --- in this case, including
the information that some features were discarded because they appear
weakly related to the response variable.  This information can only be
utilized when using a model for the joint distribution of the response
and the features, even though we are interested only in the
conditional distribution of the response given the features.

We applied this method to naive Bayes models with binary features that
are assumed to be independent conditional on the value of the binary
response (class) variable.  With these models, we can efficiently
compute the adjustment factor needed to correct for selection bias.
Crucially, we need only compute the probability that a single feature
will exhibit low correlation with the response, and then raise this
probability to the number of discarded features.  When a large number
of features are discarded, the time needed to compute the adjustment
factor for bias correction is much less that the time that would have
been needed to actually use these features.  Substantial computation
time can therefore be saved by discarding features that appear to have
little relationship with the response.

Our general method can be applied to other models and other feature
selection criteria, provided that the adjustment factor can be
computed.  Reasonably efficient computation may be possible when the
features for a case are independent given the values for a set of
latent variables, since the adjustment factor can then again be found
by raising the probability that a single feature will be discarded to
the number of features that were discarded.  However, since the values
for latent variables will not be known, the computations are more
difficult than for the naive Bayes model.  Markov chain Monte Carlo
methods will generally be needed to sample for the values of the
latent variables.  (They may be required in any case for models more
complex than the binary naive Bayes model considered in this paper.)

We have implemented such bias correction methods for two-component
mixture models of binary data, and for factor analysis models, in
which the features and the response are real valued.  The required
computations are feasible, but slower and more complex than for the
naive Bayes model.  We will report the details of these methods and
their performance in follow-on papers.  The practical utility of the
bias correction method we describe would be much improved if methods
for more efficiently computing the required adjustment factor could be
found, which could be applied to a wide class of models.

\section*{Acknowledgements}\vspace*{-8pt}
This research was supported by Natural Sciences and Engineering Research 
Council of Canada. Radford Neal holds a Canada Research Chair in 
Statistics and Machine Learning.

\section*{References}\vspace*{-8pt}

\leftmargini 0.2in
\labelsep 0in

\begin{description}
\itemsep 2pt
\item[]
  Alon, U., Barkai, N., Notterman, D.~A., Gish, K., Ybarra, S., Mack, D.,
  and Levine, A.~J.\ (1999)
  ``Broad patterns of gene expression revealed by clustering analysis
  of tumor and normal colon tissues probed by oligonucleotide arrays",
  \textit{Proceedings of the National Academy of Sciences (USA)}, vol.~96,
  pp. 6745-6750.

\item[]
  Dawid, A.~P.\  (1982) ``The well-calibrated Bayesian'', \textit{Journal
  of the American Statistical Association}, vol.~77, no.~379, pp.~605-610.

\item[]
  Guyon, I., Gunn, S., Nikravesh, M., and Zadeh, L.~A.\ (2006) \textit{Feature 
  Extraction: Foundations and Applications} (edited volume), Studies in 
  Fuzziness and Soft Computing, Volume 207, Springer.

\end{description}

\end{document}

%% file: tab.tex
\begin{tabular}
{@{}c|@{}rcc|@{}rcc|rcc|rcc@{}}
         & \multicolumn{6}{c|}{\bf  1 feature selected out of 10000 } & \multicolumn{6}{c}{\bf  10 features selected out of 10000 } \\[5pt] 
         & \multicolumn{3}{c}{Corrected} & \multicolumn{3}{c|}{Uncorrected} & \multicolumn{3}{c}{Corrected} & \multicolumn{3}{c}{Uncorrected} \\[5pt] 
Category & \mbox{~~~}\#\mbox{~} & \small{Pred} & \small{Actual} & \mbox{~~~}\#\mbox{~} & \small{Pred} & \small{Actual} & \mbox{~~}\#\mbox{~} & \small{Pred} & \small{Actual} & \mbox{~~}\#\mbox{~} & \small{Pred} & \small{Actual} \\[5pt] 
0.0 - 0.1 &    0 &    -- &    -- &    0 &    -- &    -- &    0 &    -- &    -- &  237 & 0.046 & 0.312 \\ 
0.1 - 0.2 &    0 &    -- &    -- &    0 &    -- &    -- &    3 & 0.174 & 0.000 &  349 & 0.149 & 0.444 \\ 
0.2 - 0.3 &    0 &    -- &    -- &    0 &    -- &    -- &  126 & 0.270 & 0.294 &   68 & 0.249 & 0.500 \\ 
0.3 - 0.4 &    0 &    -- &    -- & 1346 & 0.384 & 0.461 &  467 & 0.360 & 0.420 &  300 & 0.360 & 0.443 \\ 
0.4 - 0.5 & 1346 & 0.446 & 0.461 &    0 &    -- &    -- &  566 & 0.462 & 0.461 &  189 & 0.443 & 0.487 \\ 
0.5 - 0.6 &    0 &    -- &    -- &    0 &    -- &    -- &  461 & 0.554 & 0.566 &   48 & 0.546 & 0.417 \\ 
0.6 - 0.7 &  654 & 0.611 & 0.581 &    0 &    -- &    -- &  276 & 0.643 & 0.616 &  238 & 0.650 & 0.588 \\ 
0.7 - 0.8 &    0 &    -- &    -- &  654 & 0.736 & 0.581 &   97 & 0.733 & 0.742 &  180 & 0.737 & 0.567 \\ 
0.8 - 0.9 &    0 &    -- &    -- &    0 &    -- &    -- &    4 & 0.825 & 0.750 &  192 & 0.864 & 0.609 \\ 
0.9 - 1.0 &    0 &    -- &    -- &    0 &    -- &    -- &    0 &    -- &    -- &  199 & 0.943 & 0.668 \\ 
\end{tabular}

\vspace*{0.45in}
\begin{tabular}
{@{}c|@{}rcc|@{}rcc|rcc|rcc@{}}
         & \multicolumn{6}{c|}{\bf  100 features selected out of 10000 } & \multicolumn{6}{c}{\bf  1000 features selected out of 10000 } \\[5pt] 
         & \multicolumn{3}{c}{Corrected} & \multicolumn{3}{c|}{Uncorrected} & \multicolumn{3}{c}{Corrected} & \multicolumn{3}{c}{Uncorrected} \\[5pt] 
Category & \mbox{~~~}\#\mbox{~} & \small{Pred} & \small{Actual} & \mbox{~~~}\#\mbox{~} & \small{Pred} & \small{Actual} & \mbox{~~}\#\mbox{~} & \small{Pred} & \small{Actual} & \mbox{~~}\#\mbox{~} & \small{Pred} & \small{Actual} \\[5pt] 
0.0 - 0.1 &  155 & 0.067 & 0.077 &  717 & 0.017 & 0.199 &  774 & 0.018 & 0.027 &  954 & 0.004 & 0.066 \\ 
0.1 - 0.2 &  247 & 0.151 & 0.162 &  133 & 0.150 & 0.391 &   97 & 0.143 & 0.165 &   28 & 0.149 & 0.500 \\ 
0.2 - 0.3 &  220 & 0.247 & 0.286 &   70 & 0.251 & 0.429 &   63 & 0.243 & 0.302 &   13 & 0.248 & 0.846 \\ 
0.3 - 0.4 &  225 & 0.352 & 0.356 &   68 & 0.351 & 0.515 &   48 & 0.346 & 0.438 &   17 & 0.349 & 0.412 \\ 
0.4 - 0.5 &  237 & 0.450 & 0.494 &   58 & 0.451 & 0.500 &   45 & 0.446 & 0.600 &   14 & 0.449 & 0.786 \\ 
0.5 - 0.6 &  227 & 0.545 & 0.586 &   78 & 0.552 & 0.603 &   44 & 0.547 & 0.614 &   16 & 0.546 & 0.375 \\ 
0.6 - 0.7 &  202 & 0.650 & 0.728 &   77 & 0.654 & 0.532 &   53 & 0.647 & 0.698 &   16 & 0.667 & 0.812 \\ 
0.7 - 0.8 &  214 & 0.749 & 0.785 &   80 & 0.746 & 0.662 &   81 & 0.755 & 0.815 &   22 & 0.751 & 0.636 \\ 
0.8 - 0.9 &  182 & 0.847 & 0.857 &   98 & 0.852 & 0.633 &  124 & 0.854 & 0.863 &   25 & 0.865 & 0.560 \\ 
0.9 - 1.0 &   91 & 0.935 & 0.923 &  621 & 0.979 & 0.818 &  671 & 0.977 & 0.982 &  895 & 0.995 & 0.946 \\ 
\end{tabular}

\vspace*{0.45in}
\begin{tabular}[t]
{@{}c|@{}rcc@{}}
         & \multicolumn{3}{c}{\bf Complete data} \\[5pt] 
Category & \mbox{~~~}\#\mbox{~} & \small{Pred} & \small{Actual} \\[5pt] 
0.0 - 0.1 &  964 & 0.004 & 0.006 \\ 
0.1 - 0.2 &   21 & 0.145 & 0.238 \\ 
0.2 - 0.3 &    8 & 0.246 & 0.375 \\ 
0.3 - 0.4 &   10 & 0.342 & 0.300 \\ 
0.4 - 0.5 &   12 & 0.436 & 0.500 \\ 
0.5 - 0.6 &    7 & 0.544 & 1.000 \\ 
0.6 - 0.7 &   20 & 0.656 & 1.000 \\ 
0.7 - 0.8 &   13 & 0.743 & 0.846 \\ 
0.8 - 0.9 &   22 & 0.851 & 0.818 \\ 
0.9 - 1.0 &  923 & 0.994 & 0.998 \\ 
\end{tabular}